\documentstyle[amssymb,amsfonts,psfig]{amsart}

\font \roman = cmr10 at 10 true pt

\def\be#1{\begin{equation}\label{#1}}
\def\bas{\begin{align*}}
\def\eas{\end{align*}}
\def\bi{\begin{itemize}}
\def\ei{\end{itemize}}

\def\dim{{\hbox{\roman dim}}}

\newenvironment{proof}{\noindent {\bf Proof} }{\endprf\par}
\def \endprf{\hfill  {\vrule height6pt width6pt depth0pt}\medskip}
\def\emph#1{{\it #1}}
\def\textbf#1{{\bf #1}}
\theoremstyle{plain}
  \newtheorem{theorem}[subsection]{Theorem}

  \newtheorem{proposition}[subsection]{Proposition}
  
  \newtheorem{lemma}[subsection]{Lemma}
  \newtheorem{corollary}[subsection]{Corollary}
  \newtheorem{Distance conjecture}[subsection]{Distance Conjecture}
  \newtheorem{Bilinear distance conjecture}[subsection]{Bilinear Distance Conjecture}
  \newtheorem{Furstenburg problem}[subsection]{Furstenburg problem}
  \newtheorem{Discretized Furstenburg conjecture}[subsection]{Discretized Furstenburg 
  Conjecture}
  \newtheorem{Ring problem}[subsection]{Ring Problem}
  \newtheorem{Ring conjecture}[subsection]{Ring Conjecture}
  \newtheorem{Main theorem}[subsection]{Main Theorem}
  \newtheorem{Cauchy-Schwartz}[subsection]{Cauchy-Schwartz}

\theoremstyle{remark}

\theoremstyle{definition}

\include{psfig}

\begin{document}

\title[Caffarelli-Kohn-Nirenberg]{A cheap Caffarelli-Kohn-Nirenberg inequality for
Navier-Stokes equations with hyper-dissipation}

\author{Nets Hawk Katz}
\address{Department of Mathematics, Washington University 
St. Louis 63130}
\email{nets@@math.wustl.edu}

\author{Nata\v{s}a Pavlovi\'{c}}
\address{Department of Mathematics,Statistics, and
Computer Science, Chicago Il 60607-7045}
\email{natasa@@math.uic.edu}

\subjclass{???}

\begin{abstract} We prove that for the Navier Stokes
equation with dissipation $(-\Delta)^{\alpha}$, where $1<\alpha<\frac{5}{4}$,
and smooth initial data, the Hausdorff dimension of the singular
set at time of first blow up is at most $5-4\alpha$. This
unifies two directions from which one might approach the
Clay prize problem.
\end{abstract}

\maketitle

\section{Introduction}

The Navier Stokes equation with dissipation $(-\Delta)^{\alpha}$
is given by
\be{NS}{\partial u \over \partial t} + u \cdot \nabla u + \nabla p
= -(-\Delta)^{\alpha} u, \end{equation}
where $u$ is a time-dependent divergence free vector field in ${\Bbb R}^3$.
One sets the initial condition
\be{Init} u(x,0)=u_0(x) \end{equation}
where $u_0(x) \in C^{\infty}_c({\Bbb R}^3)$.

It is easy to see that classical solutions to this equation on a time interval
$[0,T]$ satisfy conservation of energy, namely that
$$||u(.,T)||_{L^2}^2 = ||u_0||_{L^2}^2 - 
\int_{0}^T \langle (-\Delta)^{\alpha} u,u \rangle.$$
The second term on the right is called the dissipation term.

Caffarelli, Kohn, and Nirenberg \cite{CKN} showed that when $\alpha=1$,
the singular set of a generalized weak solution
(a weak solution satisfying a certain generalization
of conservation of energy)  to the
system \eqref{NS},\eqref{Init} has parabolic Hausdorff dimension
at most 1. This could be considered a first step towards
showing global strong solvability. Indeed any improvement
in this upper bound on the dimension would be genuine
progress towards solution of the global solvability problem.

Another well known fact is the following. (We learned the proof
which appears in the following section from Diego Cordoba.
It appears to exist in the folklore and we have not succeeded in
attributing it. Sinai and Mattingly \cite{S&M} recently gave a
different proof.)

\begin{proposition}\label{trivial} If $\alpha>{5 \over 4}$,
one has global strong solvability for the system
\eqref{NS},\eqref{Init}. \end{proposition}

Indeed proposition \ref{trivial} could be viewed as a first
step towards the solution of the global strong solvability
for $\alpha=1$ and any improvement in the exponent ${5 \over 4}$
could be viewed as genuine progress.

The purpose of this paper is to interpolate the two results so that
these two paths to progress are unified. We prove
\begin{theorem}\label{main} If $T$ is the time of first breakdown
for the system \eqref{NS},\eqref{Init}, with 
${5 \over 4} > \alpha > 1$ then the 
Hausdorff dimension of the singular set at time $T$ is at most
$5-4\alpha$. \end{theorem}

Heuristically, the theorem ought to be a mild generalization of Proposition
\ref{trivial}. We think about the microlocal analysis of a solution $u$
in terms of coefficients $u_Q$ associated to cubes $Q$. Very roughly speaking,
the coefficient $u_Q$ should be viewed as a generalized wavelet coefficient.
The proof of Proposition \ref{trivial} goes wrong for $\alpha < {5 \over 4}$
only because of a small number of cubes with large coefficients. At any time,
the set of points contained in arbitrarily small such cubes has dimension
at most $5-4\alpha$. The difficulty lies in preventing the energy from flowing
between far away cubes near the time of blow up. This is not quite so
easy as it sounds because energy can flow from small cubes to large ones
containing them and this is a quite non-local effect.

We were also unable to directly generalize the proof of Caffarelli, Kohn, and 
Nirenberg. They rely on what they call the ``generalized energy inequality"
which is built on an amusing property of the divergence free heat equation.
Let $\phi$ be any compactly bump function in space and time and $u$ a divergence
free vector field, then
$$\int \langle ({\partial \over \partial t} - \Delta) u,\phi u \rangle dt
=\int ({1 \over 2} {\partial \over \partial t} (\langle u, \phi u \rangle)
+ \langle \nabla u, \phi \nabla u \rangle -
\langle ({{1 \over 2}{\partial \phi \over \partial t}}
+ \Delta \phi )  u , u \rangle) dt.$$
The first term represents a change in local energy. The second term represents
local dissipation. The third term is an error which can be made insignificant
by choosing $\phi$ to satisfy the backwards heat equation. It is this
method of controlling the error which we are unable to generalize. The later
proof of Lin \cite{Lin} also uses the generalized energy inequality.

To deal with localization we combine the theory of paramultiplication
following the notes of Tao \cite{Tao}. These ideas in the
use of Littlewood-Paley theory were developed by Bony (see e.g. \cite{Bony}),
and Coifman and Meyer (see e.g.\cite{Meyer}). A good exposition can  be found
in \cite{Taylor}. We combine this with the theory of type $(1,1-\epsilon)$
pseudodifferential operators which is also described in \cite{Taylor}.
The latter theory allows us to localize on cubes describing our ``wavelet
coefficients". Mimicking Caffarelli-Kohn-Nirenberg's Proposition 2 we
try to prove regularity at  points not contained in any squares in which too 
much dissipation occurs. For any such point, we derive a certain level
of regularity which we refer to as critical regularity. We find a square
centered at such a point on which we can apply a certain barrier estimate
which guarantees arbitrary regularity on its interior. This barrier
estimate can be thought of as a localized version of global existence
with small data. (One can compare this estimate with results in \cite{Tataru} and
\cite{Cannone} in which local well posedness is established with
hypotheses slightly more relaxed than ours.)
But in order for this to work,
 it is important that there not be too many small dissipating
squares on the boundary of our square. This combinatorial
issue is an ingredient which seems not
to have appeared before and which restricts us to the case $\alpha > 1$.

In section 2, we describe the proof of Proposition \ref{trivial}. In
section 3, we describe a way of measuring the Hausdorff dimension of sets.
In section 4, we introduce a dyadic model for Navier Stokes which illustrates
the ideas of our proof.
In section 5, we recall some Littlewood-Paley theory. In section 6, we describe
what happens in the way of localization of dissipation and energy
flow. Most of the general microlocal analysis occurs in this section. In
section 7, we prove a certain minimal, what we call, critical regularity
which is maintained away from a $5-4 \alpha$ dimensional set. In section
8, we describe a frequently occurring class of squares which have no highly
dissipative square on their boundary. In section 9, we apply the results
of sections 7 and 8 to bootstrap the critical regularity into arbitrarily
good regularity. 

{\bf Remark on constants and notation}

For us $\langle,\rangle$ is alway an $L^2$ pairing in space. 

Throughout this paper the expression $A \lesssim B$ means $A \leq CB$ where
$C$ is a constant. But what is a constant? Constants may depend on $T$, the
norms of the initial values, and on an $\epsilon$ which we keep fixed
throughout this paper. They may not depend on a particular square or 
a particular scale which we are estimating at.
(However sometimes we choose constants so that
there is a particular initial scale on which constants can depend.
But in that case, the point is that we make estimates depending on those 
constants  for all scales). Sometimes our hypotheses or our definitions have
constants in them. 
When this is the case, our conclusions are understood
to be a multi-parameter family of conclusions depending also on the
constants in our definitions. We will invoke propositions and lemmas using
convenient constants in the hypotheses and we do not make any attempt to
keep track of these. If we did, the paper might become unreadable.

{\bf Acknowledgements}

We would like to thank Diego Cordoba, Susan Friedlander, and Christoph 
Thiele for helpful discussions. We would especially like to thank
Terry Tao for providing us support to be at UCLA for Winter 2001 and
for the thought provoking course which he gave while we were there.
The first author was supported by NSF grant number DMS9801410.

\section{Preliminaries}

We begin by recalling the proof of Proposition \ref{trivial}.

\begin{proof}[Proof of Proposition \ref{trivial}]

We will use the standard notation that $H^{\beta}$ denotes
the $L^2$ Sobolev space over ${\Bbb R}^3$ with $\beta$ derivatives.
We will let $W^{\beta,p}$ denote the 
$L^p$ Sobolev space over ${\Bbb R}^3$ with $\beta$ derivatives.

We first recall the energy inequality obtained by pairing the
equation \eqref{NS} with $u$
\be{energy} 
{1 \over 2} {\partial (||u||_{L^2}^2) \over \partial t}
= -|| (-\Delta)^{{\alpha \over 2}} u||_{L^2}^2 
\leq -||u||_{H^{\alpha}}^2 +||u||_{L^2}^2
\end{equation}
From this, we obtain by integrating over time
(and observing that the $L^2$ norm is always positive), that if the
solution $u$ remains smooth up to time $T$, we have the estimate
\be{easyestimate}
\int_0^T ||u||_{H^{\alpha}}^2 dt \lesssim (1+T).
\end{equation}

We now pair \eqref{NS} with $(-\Delta)^{\beta}u$ in order to
estimate ${\partial(||u||_{H^{\beta}}^2) \over \partial t}$.
We obtain
\be{cleanpairing}{1 \over 2} {\partial(||(-\Delta)^{{\beta \over 2}}
u||_{L^2}^2) \over \partial t} + \langle u \cdot \nabla u,
(-\Delta)^{\beta} u \rangle = 
-||(-\Delta)^{{\alpha + \beta \over 2}}u||_{L^2}^2.
\end{equation}

Clearly, we must estimate the nonlinear problem term
$$ \langle u \cdot \nabla u,(-\Delta)^{\beta} u \rangle 
=\langle (-\Delta)^{{\beta \over 2}}( u \cdot \nabla u),
(-\Delta)^{{\beta \over 2}} u \rangle.$$
Since $u$ is divergence free, we can bound 
the absolute value of the above by
$$ ||u||_{W^{1,p}} ||u||_{W^{\beta,q}} ||u||_{H^{\beta}} ,$$
with ${1 \over p} + {1 \over q} ={1 \over 2}$. Now since
we are assuming $\alpha > {5 \over 4}$, we obtain by the
Sobolev imbedding theorem,
$$ ||u||_{W^{1,p}} ||u||_{W^{\beta,q}} ||u||_{H^{\beta}} 
\lesssim ||u||_{H^{\alpha}} ||u||_{H^{\alpha+\beta}}
||u||_{H^{\beta}},$$
since we must spend $> {3 \over 2}$ derivatives to get all
three terms in $L^2$ and we have spent $2 \alpha-1$.

Applying Cauchy-Schwartz, we get immediately
$$||u||_{H^{\alpha}} ||u||_{H^{\alpha+\beta}}
||u||_{H^{\beta}} \leq \delta ||u||_{H^{\alpha+\beta}}^2
+ {1 \over \delta}  ||u||_{H^{\alpha}}^2 ||u||_{H^{\beta}}^2.$$

Combining this with \eqref{cleanpairing}, we get
$$ {\partial ||u||_{H_{\beta}}^2 \over \partial t}
\lesssim ||u||_{H^{\alpha}}^2 ||u||_{H^{\beta}}^2 + ||u||_{L^2}^2.$$

In turn, combining this with \eqref{easyestimate} and with 
Gronwall's inequality gives global solvability. \end{proof}

\section{Hausdorff dimension}

The purpose of this section is to recall the definition of Hausdorff
dimension and to explain a simple construction which automatically
produces sets below a certain desired dimension.

Given any set $E \subset {\Bbb R}^n$, we define its $d$-dimensional
Hausdorff measure ${\cal H}^d(E)$ by 
$$ {\cal H}^d(E) =
\lim_{\delta \longrightarrow 0} 
\inf_{C \in {\cal C}_{\delta}(E)} \sum_{B \in C} r(B)^d.$$
Here ${\cal C}_{\delta}$ is the set of all coverings of $E$
by balls of radius less than or equal to $\delta$, and if $B$
is a ball, $r(B)$ is its radius.

We say the Hausdorff dimension of $E$ is the infimum of the
set of all $d$ for which ${\cal H}^d(E)$ is 0. Note
that for the purposes of the definition of dimension
it doesn't matter what basis of objects we cover by 
(balls, squares, ellipsoids of a fixed multiplicity)
as long as they are a regular family. (See \cite{Stein}.) 

All of our arguments are built on the notion of
discretizing at a certain scale $2^{-j}$. We would 
like to be able to think of sets of a certain
Hausdorff dimension as sets of $2^{jd}$ elements
at scale $2^{-j}$. We make this intuition precise
for our purposes by the following.

\begin{lemma}\label{limsup} Let $A_1, \dots, A_j, \dots$
be a sequence of collections of balls
in ${\Bbb R}^n$ so that each element of $A_j$ has radius $2^{-j}$.
Suppose that $\#(A_j) \leq 2^{jd}$. Define
$$E=\limsup_{j \longrightarrow \infty} A_j,$$
to be the set of points in infinitely many of 
the  $\cup_{B \in A_j} B$'s. Then
the Hausdorff dimension of $E$ is at most $d$.
\end{lemma}

\begin{proof}
Since the Hausdorff dimension of E is the infimum of the set of all
$\gamma$ for which ${\cal H}^{\gamma}(E) = 0$, it is enough to prove that
${\cal H}^{\gamma}(E) = 0$ for all $\gamma > d$.

Pick $j$ such that $2^{-j} < \delta$. We can cover  E by the $\cup_{k>j}
\cup_{B \in A_k} B$.

Now
$$ {\cal H}^{\gamma}(E) \leq \sum_{k>j} 2^{kd}(2^{-k})^{\gamma},$$
which limits to zero as $j$ goes to $\infty$ whenever $\gamma > d$. \end{proof}

\section{Dyadic heuristic}

Before going into the details of the proof of Theorem \ref{main}, we present a
heuristic to illustrate the main ideas. 

We shall replace our vector valued function $u$ by a scalar valued one.
We shall represent it by a wavelet expansion. This shall be left 
rather abstract. We do not care which system of wavelets we use
precisely so long as the size of the coefficients of a function
with respect to the basis accurately reflect the regularity properties
of the function. We denote
our orthonormal family of wavelets $\{ w_Q \}$, 
with $ w_Q $ the wavelet associated to the spatial dyadic cube $Q$ of
sidelength $2^{-j(Q)}$. Then $u$ can be represented as:
$$ u = \sum_Q u_Q w_Q, $$
where $u_Q$ is the wavelet coefficient associated to the cube $Q$. 

Until the end of this section we reserve notation $w_Q$ to denote the wavelet
localized on a cube $Q$, while $u_Q$, and $v_Q$ are reserved for wavelet
coefficients associated to a cube $Q$. 

We need a few definitions in order to make our model precise.
A cube $Q$ in ${\Bbb R}^{3}$ is called a dyadic cube if its sidelength is
integer power of 2, $2^{l}$, and the corners of the cube are on the
lattice $2^{l} {\Bbb Z}^{3}$. A dyadic child of a dyadic cube $Q$
with sidelength $2^{l}$ is any of $2^{3}$ dyadic cubes $Q^{\prime}$
contained in cube $Q$ and such that their sidelengths are $2^{l-1}$.

We will call $(Q, Q^{\prime}, Q^{\prime \prime})$ a cascade if
$Q^{\prime}$ is a dyadic child of $Q$, and $Q^{\prime \prime}$ is a
dyadic child of $Q^{\prime}$.

Now we are ready to define an operator which will mimick the behavior of
nonlinear term $u \cdot \nabla u$. Note that we
have 
\begin{equation}
||w_Q||_{L^{\infty}} \sim 2^{\frac{3j(Q)}{2}}. \label{wav-infty} 
\end{equation}

On the other hand 
\begin{equation}
||\nabla w_Q||_{L^{2}} \sim 2^{j}, \label{wav-nabla}
\end{equation}
since our wavelets $w_Q$'s are orthonormal and localized to the frequencies
around $2^{j}$. 

Having in mind (\ref{wav-infty}) and (\ref{wav-nabla})
we define a cascade down operator as follows:
$$ c_d (u,v) = \sum_{cascades} 2^{\frac{3j(Q)}{2}} 2^{j(Q)} u_Q v_{Q^{\prime}}
w_{Q^{\prime \prime}}. $$ 

Similarly we define a cascade up operator as:
$$ c_u (u,v) = \sum_{cascades} 2^{\frac{3j(Q)}{2}} 2^{j(Q)} u_{Q^{\prime
\prime}} v_{Q^{\prime}} w_Q. $$ 

(Note that the choice of definition of the cascade down operator is somewhat
arbitrary. We are happy to take any bilinear operator with few coefficients
which shifts energy to high frequency only a few scales at a time. Once
cascade down is chosen, cascade up must be one of its adjoints.)

Also let us define 
$$ c(u,v) =  c_d (u,v) -  c_u (u,v). $$

Obviously, 
$$ \langle c_d (u,u), u \rangle = \langle c_u (u,u), u \rangle, $$
which implies 
\begin{eqnarray}
\langle c(u,u), u \rangle = 0. \label{cons}
\end{eqnarray}

Having defined operator $c(u,v)$ we can speak about dyadic version of
Euler as well as Navier-Stokes equations. More precisely by dyadic Euler
equation we mean:
$$ \frac{d u}{d t} + c(u,u) = 0.$$

For the remainder of this section, we define
 Laplacian as  $ \Delta (w_Q) = 2^{2j} w_Q $. Now we introduce
dyadic version of Navier-Stokes equation as:
$$ \frac{d u}{d t} + c(u,u) + \Delta u = 0.$$

Also we could speak about dyadic Navier-Stokes equation with
hyper-dissipation by which we mean:
\begin{eqnarray}
\frac{d u}{d t} + c(u,u) + (\Delta)^{\alpha} u = 0. \label{dyadichyp}
\end{eqnarray}

A simple consequence of (\ref{cons}) is conservation of energy for all
three equations, which is an important feature of these equations
preserved in the dyadic model described here.

Now we consider dyadic model for the Navier-Stokes equation with
hyper-dissipation (\ref{dyadichyp}). We would like to estimate Hausdorff
dimension of the set of singular points at the first time of blow up, 
$T$. The nonlinear term $c(u,u)$ on scale $j$ looks like $2^{\frac{5j}{2}}
u_{Q}^{2}$ (if we imagine for a moment that all neighboring cubes have
coefficients of roughly the same size), while 
dissipation term gives decay like $2^{2 \alpha j}
u_Q$. This means that as long as $u_Q < 2^{-\frac{j}{2} (5-4\alpha)}$, the
growth of $u_Q$ is under control. But let us check what happens if $u_Q >
2^{-\frac{j}{2} (5-4\alpha)}$. (This idea is made more precise by
Lemma \ref{Barrier-estimate}. The reader is encouraged to work out its
simpler dyadic analogue). We call this bound on the coefficients
critical regularity and this is what we shall prove away from a
$5-4\alpha$ dimensional set.

We can rewrite equation (\ref{dyadichyp}) in terms of wavelets
coefficients as follows:
\begin{eqnarray}
\frac{d u_{Q}}{d t} =  \sum_{(Q^{\prime},Q^{\prime \prime})}
c_{(Q,Q^{\prime} ,Q^{\prime \prime})}
2^{\frac{5j(Q)}{2}} u_{Q^{\prime}} u_{Q^{\prime
\prime}} - 2^{2 \alpha j(Q)} u_Q , \label{coefNS}
\end{eqnarray}

where $c_{(Q,Q^{\prime},Q^{\prime \prime})}=-1$ if $(Q,Q^{\prime},Q^{\prime \prime})$
is a cascade, where $c_{(Q,Q^{\prime},Q^{\prime \prime})}=32$ if
$(Q^{\prime}, Q^{\prime \prime},Q)$ is a cascade and vanishes otherwise.

Having assumed $u_Q \gtrsim 2^{-\frac{j}{2} (5-4\alpha)}$ for some time $t$
and assuming that at the initial time $t=0$ it is much smaller, by
the smoothness assumption on the initial condition, 
we integrate
(\ref{coefNS}) in time on the interval $[0,T]$ and obtain
for one of the choices of $(Q^{\prime} ,Q^{\prime \prime})$ giving a non-vanishing
coefficient:
$$ 2^{\frac{5j}{2}} \int_{0}^{T} |u_{Q^{\prime}} u_{Q^{\prime \prime}}| dt \gtrsim
2^{-\frac{j}{2} (5-4\alpha)}, $$ 
which by Cauchy-Schwartz implies:
$$ 2^{\frac{5j}{2}} (\int_{0}^{T} u_{Q^{\prime}}^{2} dt)^{1 \over 2}
(\int_{0}^{T} u_{Q^{\prime \prime}}^{2} dt)^{1 \over 2}  \gtrsim
2^{-\frac{j}{2} (5-4\alpha)}. $$

The last expression can be rewritten as:
$$2^{2j \alpha} (\int_{0}^{T} u_{Q^{\prime}}^{2} dt)^{1 \over 2}
(\int_{0}^{T} u_{Q^{\prime \prime}}^{2} dt)^{1 \over 2}  \gtrsim
2^{-j(5-4\alpha)}. $$

This can happen if either 
\begin{eqnarray}
2^{2j \alpha} (\int_{0}^{T} u_{Q^{\prime}}^{2} dt)^{1 \over 2} \gtrsim
2^{-j(5-4\alpha)},
\label{possib1}
\end{eqnarray}
or
\begin{eqnarray}
2^{2j \alpha} (\int_{0}^{T} u_{Q^{\prime \prime}}^{2} dt)^{1 \over 2} \gtrsim
2^{-j(5-4\alpha)}.
\label{possib2}
\end{eqnarray}

However having in mind conservation of energy we have
$$ 2^{2j \alpha} \int_{0}^{T} \sum_{Q\; at \; scale \;j} u_{Q}^{2} dt 
\lesssim 1.$$
Thus we conclude that (\ref{possib1}) or  (\ref{possib2})
could happen in at most $\lesssim 2^{j(5-4\alpha)}$ cubes $Q$. Now we invoke Lemma
\ref{limsup} and conclude that the Hausdorff dimension of the set of
points of the equation (\ref{dyadichyp})
at which critical regularity fails is at most $5-4\alpha$.

In the remainder of the paper, we shall prove the same result for
the actual Navier-Stokes with hyper-dissipation and we shall prove that
this critical regularity implies stronger regularity at points away
from the bad set.

We would like to remark that the dyadic model seems quite interesting in that 
it possesses all the features of Navier Stokes which we shall use in this
paper. But moreover the nonlinear term posesses a built-in dispersive
feature - namely that cascade down passes energy 
equally to all the squares $Q^{\prime \prime}$
which are children of $Q^{\prime}$ and grandchildren of $Q$. This dispersive
feature might be useful in proving global solvability for the model. We believe
that working out such an argument would be a good first step towards the
Clay problem.

\section{Littlewood Paley theory and Pseudodifferential operators}

We shall use a standard Littlewood Paley partition of
frequency space. That is, we shall have Fourier multipliers
$P_j$ (on $L^2({\Bbb R}^3)$)   so that their symbols $p_j(\xi)$ are smooth and supported
in ${2 \over 3} 2^j < |\xi| < 3 (2^j)$, so that $p_j(\xi)=p_0(2^{-j} x)$
 and so that
$$\sum_j p_j(\xi)=1.$$

We define $\tilde P_j=\sum_{k=-2}^2 P_{j+k}$. Notice that 
$\tilde P_j P_j=P_j$, since $\tilde P_j$ is the sum of all
Littlewood-Paley projections the support of whose symbols intersects
the support of $p_j(\xi)$. We can analogously define the symbols
$\tilde p_j=\sum_{k=-2}^2 p_{j+k}$. 

The Littlewood-Paley operators satisfy the inequality

\begin{proposition}\label{Bernstein} 
$$||P_j f||_{L^{\infty}} \lesssim 2^{3j \over 2} ||P_j f||_{L^2}.$$
\end{proposition}

\begin{proof}
Let $\phi_j = \check ( \tilde p_j)$. Then (since we are working in ${\Bbb R}^3$) we have
$$\phi_j(x)=2^{3j} \phi_0(2^j x).$$
Now $\tilde P_j f = \phi_j * f$ so by Young's inequality, we have
$$||P_j f||_{L^{\infty}} = ||\tilde P_j P_j f||_{L^{\infty}} \leq
||\phi_j||_{L^2} ||P_j f||_{L^2},$$
which is finite since $\phi_j$ is a Schwartz function. But by the
relation between $\phi_j$ and $\phi_0$, we have
$$||\phi_j||_{L^2} = 2^{{3j \over 2}} ||\phi_0||_{L^2},$$
which proves the proposition.
\end{proof}

The idea of Littlewood Paley theory is that $P_j f$ is like
a combination of wavelets supported on cubes of length
$2^{-j}$ and that the previous proposition is sharp only
when a large proportion of the $L^2$ energy of $P_j f$ is
concentrated in a single one of these cubes. Hence one is led
to try to localize $P_j f$ in space on such cubes and this
is marginally consistent with the Heisenberg uncertainty
principle. However we have some room for error in our dimension
estimate since it is a closed condition. Therefore
we fix an $\epsilon > 0$ and never try to localize better
than within $2^{-j(1-\epsilon)}$.

Our localization is nearly perfect
if we neglect negligible quantities. For the 
remainder of the paper, whenever we are at
scale $j$, we neglect quantities of size 
$\lesssim 2^{-100j}$ since they will not affect our estimates.
Similarly we neglect operators whose norms are smaller than
$2^{-100 j}$ provided they will only be applied to functions
whose norms are $\lesssim 1$.
The choice of 100 is arbitrary. It is large enough so that it
does not affect the $L^{\infty}$ norm of $u$. But in fact
our techniques for showing quantities are negligible rely on Schwartz
function properties and could give arbitrary exponent with loss in
the constant. The current approach works well principally
because of conservation of energy, since the $L^2$ norm of $u$
is certainly $\lesssim 1$.

For any cube $Q$ with sidelength greater than $2^{-j(1-2\epsilon)}$,
we define a bump function $\phi_{Q,j}$ 
which is positive, is bounded above by
$1$,  equals $1$ on $Q$ and is 0 outside
of $(1+2^{-j\epsilon})Q$. Further, we require for each multiindex
$\alpha$ that there is a constant $C_{\alpha}$ independent of $Q$
so that
\be{symbests}|D^{\alpha} \phi_{Q,j}| \leq C_{\alpha} 2^{|\alpha| j (1 -\epsilon)}.
\end{equation}
We say that any bump function which satisfies the estimates
\eqref{symbests} is of type $j$.

The map $\phi_{Q,j} P_j$ acts much like a
projection, and we shall treat $||\phi_{Q,j} P_j f||_{L^2}$ 
as if it were a wavelet coefficient. When we deal with a cube
$Q$ of sidelength $2^{-j(1-2 \epsilon)},$ we shall define $j(Q)=j$
and we shall denote $||\phi_{Q,j(Q)} P_j f||_{L^2}$ by $f_Q$. Further
if $j(Q)=j$, we say that $Q$ is at level $j$.

\begin{proposition}\label{freqloc}
Given $f$ with $||f||_{L^2} \lesssim 1$, and $\phi$ a bump function of
type $j$
we have that the quantity
$$||\phi  P_j f - 
\tilde P_j \phi P_j f||_{L^2}$$
is negligible.
\end{proposition}

\begin{proof}
Let us define $\phi=\phi_1 + \phi_2,$ where
$\hat \phi_2(\xi)=\chi_{\{|\xi| > {1 \over 100} 2^j\}} \hat \phi(\xi)$.
We have 
$$\phi  P_j f - 
\tilde P_j \phi P_j f=\phi_1 P_j f - \tilde P_j \phi_1 P_j f
+\phi_2 P_j f - \tilde P_j \phi_2 P_j f.$$

By our estimates on the derivatives of $\phi$, we can get
$$||\phi_2(\xi)||_{L^{\infty}}=negligible,$$
while, because of the Fourier transform supports,
$$\tilde P_j \phi_1 P_j f =\phi_1 P_j f.$$
Thus the proposition is proved. \end{proof}

Combining 
propositions \ref{Bernstein} and \ref{freqloc}, we get the following extremely
useful lemma.

\begin{lemma}\label{superbernstein}
Let $\phi$ be a bump function of type $j$. Then
$$||\phi P_j f||_{L^{\infty}} \lesssim 2^{{3j \over 2}} ||\phi P_j f||_{L^2}
+ negligible.$$
\end{lemma}

\begin{proof} We neglect negligible terms. Then we have $\phi P_j f = 
\tilde P_j \phi P_j f$, by proposition \ref{freqloc}. 
We estimate $||\tilde P_j \phi P_j f||_{L^{\infty}}$ by proposition \ref{Bernstein}.
\end{proof}

\begin{proposition}\label{spacloc} For any cube $Q$, we have that
for any $f$ with $||f||_{H^{\beta}} \lesssim 1$ with $\beta <10$ that
$$(1 - \phi_{(1 + 2^{-\epsilon j(Q) \over 2})Q,j(Q)}) P_j \phi_{Q,j(Q)} f,$$
is negligible. \end{proposition}

\begin{proof} 
$$(1 - \phi_{(1 + 2^{-\epsilon j(Q) \over 2})Q,j(Q)}) P_j \phi_{Q,j(Q)} ,$$
is a composition of type $(1,1-\epsilon)$ pseudodifferential operators whose
symbols have disjoint support. Thus it is smoothing.
\end{proof}

Notice that the above proposition really says that we can move bump functions
across Littlewood Paley projections as long as the bump functions proliferate
and increase in support. In other words, the proposition can be rewritten as
$$\phi_{(1 + 2^{-\epsilon j(Q) \over 2})Q,j(Q)} P_j \phi_{Q,j(Q)} =
P_j \phi_{Q,j(Q)} + negligible.$$
Similarly, this can be used to remove inconvenient bump functions, provided there is
a smaller bump function already present in the expression.

\begin{proposition}\label{covering}
Let ${\cal C}$ be a covering of a set $E$ by cubes of sidelength $2^{-j(1-2\epsilon)}$
for some fixed $j$. Then for any $f \in L^2$, we have
$$||\chi_E P_j f||_{L^2}^2 \leq \sum_{Q \in \cal C} f_Q^2.$$
\end{proposition}

\begin{proof} We simply observe that
$$ \sum_{Q \in \cal C} f_Q^2=\int (\sum_{Q \in \cal C} \phi_{Q,j}^2) |P_j f|^2
\geq \int \chi_E |P_j f|^2 .$$
\end{proof}

Essentially what we have done in this section is to use 
approximations to projections which
are uniformly pseudodifferential operators of type $(1,1-\epsilon)$. 
The negligible interaction of distant squares can be just as well derived
from the asymptotic formula for composition of such operators. The composition
of two $(1,1-\epsilon)$ operators whose symbols have disjoint support is
infinitely smoothing and hence negligible.

\section{Dissipation, Trichotomy, and Forbidden squares}

In this section, with $u$ a solution to the Navier Stokes with
hyperdissipation, we write down estimates involving the
growth of the $u_Q$'s. This is the only section in which
we work directly with the equation. By the end of this section,
we will have reduced our problem to a combinatorial one about
ordinary differential equations on graphs. This section is
presented with apologies to Terry Tao, since much of it
is based on a course he gave at UCLA in Winter 2001.

Recall we have
\be{NS+}{\partial u \over \partial t} + u \cdot \nabla u + \nabla p
= -(-\Delta)^{\alpha} u, \end{equation}
together with
\be{divfree} \nabla \cdot u =0 \end{equation}

In light of \eqref{divfree}, we can rewrite \eqref{NS+} as

\be{NS++}{\partial u \over \partial t} + T(u \cdot \nabla u) 
= -(-\Delta)^{\alpha} u, \end{equation}

where $T$ is the projection into divergence free vector fields.
The operator $T$ is a singular integral operator and is also
a Fourier multiplier. We pick a cube $Q$ of side length
$2^{-j(1-2\epsilon)}$ and compute the $L^2$ pairing
of the equation with $P_j \phi_{Q,j}^2 P_j u$. We obtain
the energy estimate

\be{energy1} {1 \over 2} {d \over dt} u_Q^2=
\langle -T(u \cdot \nabla u),P_j \phi_{Q,j}^2 P_j u \rangle
-\langle (-\Delta)^{\alpha} u,P_j \phi_{Q,j}^2 P_j u  \rangle
\end{equation}

We shall estimate the two terms on the right hand side of
\eqref{energy1} separately.

As before we use the notation 
$$\tilde P_j = \sum_{k=-2}^2 P_{j+k}.$$
This has the advantage that
$$P_j = \tilde P_j P_j.$$

We make some definitions. 
We define for each cube $Q$, the set
${\cal N}^1(Q)$, the nuclear family of $Q$ to be a union of sets
$ A_Q,B_Q,C_Q,D_Q,E_Q$, where $A_Q,B_Q,C_Q,D_Q,E_Q$ are covers
of $\tilde Q=(1 + 2^{-\epsilon j \over 4}) Q$ by fewer 
than 1024 cubes each at
levels respectively $j-2,j-1,j,j+1,$ and $j+2$. We define recursively
${\cal N}^l(Q)$ to be the union of all ${\cal N}^1(Q^{\prime})$ for all
$Q^{\prime} \in {\cal N}^{l-1}(Q)$. Thus in particular, we have
$$\#({\cal N}^l(Q)) \leq 2^{13 l},$$
since $5(2^{10}) < 2^{13}$.

\begin{proposition}\label{diss} Let $Q$ be a cube and $j=j(Q)$. Then
\be{diss1}
\langle (-\Delta)^{\alpha} u,P_j \phi_{Q,j}^2 P_j u  \rangle
\geq  c 2^{2\alpha j}u_Q^2 - C 2^{(2 \alpha -2\epsilon)j}
\sum_{Q^{\prime} \in {\cal N}^{1}(Q)}
u_{Q^{\prime}}^2 - negligible  \end{equation}
\end{proposition}

\begin{proof}
Note that
$$\langle (-\Delta)^{\alpha} u,P_j \phi_{Q,j}^2 P_j u \rangle=
\langle (-\Delta)^{\alpha} \phi_{Q,j} P_j u, \phi_{Q,j} P_j u \rangle
+\langle [\phi_{Q,j} P_j,(-\Delta)^{\alpha}] u,\phi_{Q,j} P_j u \rangle
=X+Y.$$

Note that by proposition \ref{freqloc}, we have
$$X \gtrsim 2^{2\alpha j} u_Q^2 - negligible.$$

To estimate $Y$, we observe that $[\phi_{Q,j} P_j,(-\Delta)^{\alpha}]$
is of order $2\alpha-\epsilon$ since $\phi_{Q,j} P_j$ is of type
$(1,1-\epsilon)$. Further by proposition \ref{freqloc}, we have
$$[\phi_{Q,j} P_j,(-\Delta)^{\alpha}]\tilde P_j -
[\phi_{Q,j} P_j,(-\Delta)^{\alpha}] =negligible.$$
Even further, applying proposition \ref{spacloc}, we get
$$Y= \langle [\phi_{Q,j} P_j,(-\Delta)^{\alpha}]\phi_{\tilde Q,j} \tilde P_j
 u,\phi_{Q,j} P_j u \rangle + negligible.$$
Now applying the mapping properties of operators of order $2\alpha - \epsilon$
and Cauchy Schwartz, we get for any number $K$,
$$|Y| \lesssim {1 \over K} 2^{2\alpha j} u_Q^2 + 
K 2^{(2\alpha-2\epsilon)j} ||\phi_{\tilde Q,j} \tilde P_j u||_{L^2}^2.$$
Now applying proposition \ref{covering}, we get the desired result.
\end{proof}

Let us define for every $l$,
$$u_{{\cal N}^l(Q)}^2 = \sum_{Q^{\prime} \in {\cal N}^l(Q)} u_{Q^{\prime}}^2.$$

Then we have
\begin{corollary}\label{neighdiss} Let $Q$ be a cube and $j=j(Q)$. Then for any $l$,
\be{diss2} \sum_{Q^{\prime} \in {\cal N}^l (Q)}
\langle (-\Delta)^{\alpha} u,P_j \phi_{Q^{\prime},j}^2 P_j u  \rangle
\geq  c 2^{2\alpha j}u_{{\cal N}^l (Q)}^2 - C 2^{(2 \alpha -2\epsilon)j}
u_{{\cal N}^{l+1} (Q)}^2 - negligible  \end{equation}
\end{corollary}

\begin{proof} Simply sum proposition \ref{diss} over ${\cal N}^l (Q)$.
\end{proof}

We want to use Proposition \ref{diss} to fight against the
growth of $u_Q$. Thus the term $2^{(2\alpha-2\epsilon) j} 
u_{Q^{\prime}}^2$ would appear to be a serious nuisance.
However we are able
to show

\begin{lemma}\label{goodneighbours}
For any interval of time $J \subset [0,T]$ and any cube $Q$
we may find an $l < {400 \over \epsilon}$  so that 
\be{gen}\int_J u_{{\cal N}^l (Q)}^2 dt
+ negligible \gtrsim 2^{-\epsilon j}
\int_J u_{{\cal N}^{l+1} (Q)}^2 dt.\end{equation}
\end{lemma}

\begin{proof}
By conservation of energy, for every $l \leq {400 \over \epsilon}$ and every $t$
we have,
$$u_{{\cal N}^{l+1} (Q)}^2 \lesssim 1.$$
Since our constants can depend on $T$ this means
$$\int_J u_{{\cal N}^{l+1} (Q)}^2 dt \lesssim 1.$$
Now suppose the lemma is false. Applying the opposite of \eqref{gen},
for $l=1 \dots {400 \over \epsilon},$ we get
$$\int_J u_{{\cal N}^{2} (Q)}^2 dt + negligible \lesssim
2^{-399j} = negligible,$$
which implies that \eqref{gen} holds for $l=1$.
\end{proof}

We will be able to apply Lemma \ref{goodneighbours} together with Corollary
\ref{neighdiss} to show that for any cube and any interval in time, there
is an $l$th iterated nuclear family for small $l$ which is undergoing
dissipation.

Now we turn our attention to the term
$$G_Q=\langle -T(u \cdot \nabla u),P_j \phi_{Q,j}^2 P_j u \rangle.$$
We rewrite it as
$$G_Q=\langle 
-\phi_{Q_j} \tilde P_j T P_j(u \cdot \nabla u),
 \phi_{Q,j} P_j u \rangle.$$

Now we  use the ``trichotomy". We can write
$$P_j(u \cdot \nabla u) = H_{j,lh} + H_{j,hl} + H_{j,hh} + H_{loc},$$
where the low-high part is given by
$$H_{j,lh}=\sum_{k< j - {1000 \over \epsilon}}
P_j((P_k u) \cdot \tilde P_j \nabla u),$$
the high-low part is given by
$$H_{j,hl}=\sum_{k< j - {1000 \over \epsilon}}
P_j((\tilde P_j u )\cdot P_k \nabla u),$$
the high-high part is given
$$H_{j,hh}=
\sum_{k> j + {1000 \over \epsilon}}
P_j((\tilde P_k u) \cdot P_k \nabla u)
+ \sum_{k> j + {1000\over \epsilon}}
P_j(( P_k u) \cdot \tilde P_k \nabla u),
 ,$$
(technically that was the end of the trichotomy) and the local part
is given by
$$H_{loc}=
\sum_{ j- {1000 \over \epsilon}<k< j + {1000 \over \epsilon} }
P_j((\tilde P_k u) \cdot P_k \nabla u)
+ \sum_{ j- {1000 \over \epsilon}<k< j + {1000 \over \epsilon}  }
P_j(( P_k u) \cdot \tilde P_k \nabla u),
 .$$
Now we break up $G_Q$ in the obvious way:
$$G_Q =  G_{Q,lh} + G_{Q,hl} + G_{Q,hh} + G_{Q,loc},$$
where
$$G_{Q,lh}=\langle 
-\phi_{Q,j} \tilde P_j T H_{j,lh},
 \phi_{Q,j} P_j u \rangle.$$
$$G_{Q,hl}=\langle 
-\phi_{Q,j} \tilde P_j T H_{j,hl},
 \phi_{Q,j} P_j u \rangle.$$
$$G_{Q,hh}=\langle 
-\phi_{Q,j} \tilde P_j T H_{j,hh}
 \phi_{Q,j} P_j u \rangle.$$
$$G_{Q,loc}=\langle 
-\phi_{Q,j} \tilde P_j T H_{j,loc},
 \phi_{Q,j} P_j u \rangle.$$

Now fixing $Q$  and $j$ we would like to estimate each of 
$G_{Q,lh},G_{Q,hl},G_{Q,hh},$ and $G_{Q,loc}$.

We take a moment to be careful about how we localize. For any $k < 
j$, define
$Q_k=2^{(j-k)(1-2\epsilon)} Q$. For any $k \geq j$ (including $j$)
define $Q_k=(1 + 2^{-\epsilon k \over 2}) Q$. We retain this notation
only for this section and the following one.

Then we have

\begin{lemma} \label{highlowlowhigh} For any $\delta > 0$
$$|G_{Q,lh}| + |G_{Q,hl}| \lesssim \sum_{k=\delta j}^{j-{1000 \over \epsilon}} 
   2^{{3k \over 2} + j} u_{Q_k} u_{{\cal N}^1(Q)} u_Q + 
2^{j(1 + {3 \delta \over 2})} u_{{\cal N}^1(Q)} u_Q +
negligible.$$
\end{lemma}

\begin{proof}
We consider for $k < j -{1000 \over \epsilon}$, the relevant expression
$$G_{Q,lh,k}=\langle 
-\phi_{Q,j} P_j T (P_k u \tilde P_j \nabla u), \ 
 \phi_{Q,j} P_j u \rangle.$$
We divide this into two cases which are $k < \delta j$ and $k \geq \delta j$.
In the second case, we use the idea of proposition \ref{spacloc} (we just have
to use a higher degree of smoothing than is used in that proposition) to observe
$$G_{Q,lh,k}=\langle 
-\phi_{Q,j} P_j T ((\phi_{Q_k,k} P_k) u \cdot \nabla(\phi_{Q_j,j}\tilde P_j  
u)),\ \phi_{Q,j} P_j u \rangle + negligible.$$
(This is because $\nabla$ acts on $\phi_{Q_j,j}$ only where it
has a negligible effect on the whole quantity by proposition \ref{spacloc}.)
Now we simply observe by proposition \ref{superbernstein} that
$$|| \phi_{Q_k,k} P_k u ||_{L^{\infty}} \lesssim 2^{{3k \over 2}} u_{Q_k}
+ negligible,$$
and by proposition \ref{spacloc},  the proof of proposition \ref{freqloc} (to
control the action of $\nabla$) and proposition \ref{covering} that
$$ ||\nabla \phi_{Q_j,j} \tilde P_j u||_{L^2} \lesssim 2^j u_{{\cal N}^1(Q)}
+ negligible.$$

On the other hand for the first case, there is no point in localizing to
$Q_k$ because $k$ is too small and so the error is too big. Thus in in this
case, we simply estimate 
$$||P_k u ||_{L^{\infty}} \lesssim 2^{{3\delta j \over 2}} u_{Q_k}.$$
Summing these estimates gives the desired bound for $G_{Q,lh}$.
The bound for $G_{Q,hl}$ proceeds likewise (and gives a better estimate
since the derivative falls on the level $k$ term.)

\end{proof}

\begin{lemma} \label{neighbours}
$$|G_{Q,loc}| \lesssim 
2^{5j \over 2} u_Q u_{{\cal N}^{{1000 \over \epsilon}}(Q)}^2 + negligible.$$
\end{lemma}

\begin{proof}
From the definition of $G_{Q,loc}$ we have:
$$G_{Q,loc}=\sum_{l=-2}^{2} [
\sum_{k=j-{1000 \over \epsilon}}^{j+ {1000 \over \epsilon}}
\langle \phi_{Q,j} P_j T(P_{k+l} u \cdot \nabla P_{k} u), \phi_{Q,j} P_j u
\rangle
+\sum_{k=j-{1000 \over \epsilon}}^{j+ {1000 \over \epsilon}}
\langle \phi_{Q,j} P_j T(P_{k} u \cdot \nabla P_{k+l} u), \phi_{Q,j} P_j u \rangle ].$$

Since the above sum has only $\lesssim 1$ many terms by applying Proposition \ref{spacloc}, 
we observe that for some particular values of $k,l$
$$|G_{Q,loc}|\lesssim |\langle \phi_{Q,j} P_j T(\phi_{Q_j,j} P_k u \cdot 
\phi_{Q_j,j} \nabla P_{k+l} u),P_j \phi_{Q,j} u \rangle |+ negligible.$$
Now we apply proposition \ref{superbernstein} together with Cauchy Schwartz to obtain 
$$ ||\phi_{Q_j,j} P_k u||_{L^{\infty}} \lesssim 2^{{3j \over 2}}
||\phi_{Q_j,j} P_k u||_{L^2} + negligible.$$
Using proposition \ref{covering}, we observe that
$$ ||\phi_{Q_j,j} P_k u||_{L^2} \lesssim u_{{\cal N}^{{1000 \over \epsilon}}(Q)}.$$
Finally direct calculation shows
$$||\phi_{Q_j,j} \nabla P_{k+l} u||_{L^2} \lesssim ||\nabla 
\phi_{(1+2^{-2 \epsilon j \over 3})Q_j,j} P_{k+l} u|| \lesssim 
2^j u_{{\cal N}^{{1000 \over \epsilon}}(Q)} + negligible,$$
where the first inequality comes from the fact that 
$\phi_{(1+2^{-2 \epsilon j \over 3})Q_j,j}=1$ on the support of $\phi_{Q_j,j}$
and the second inequality comes again from proposition \ref{covering}.
Combining all these estimates proves the lemma.

\end{proof}

\begin{lemma} \label{highhigh}
$$|G_{Q,hh}| \lesssim 
\sum_{k > j + {1000 \over \epsilon}} u_Q 2^{{3j \over 2} + k}
||\phi_{Q_j,j} P_k u||_{L^2}^2.$$
\end{lemma}

\begin{proof}
We begin similarly to before by estimating 
$$\langle \phi_{Q,j} P_j T (P_k u \cdot \nabla \tilde P_k u), \phi_{Q,j} 
P_j u \rangle.$$
With the other terms one can proceed likewise.
Now, similarly to what we have already done, we observe that we can write this as,
$$\langle \phi_{Q,j} P_j T (\phi_{Q_j,j} P_k u \cdot \nabla
\phi_{Q_j,j} \tilde P_k u), \phi_{Q,j} P_j u \rangle + negligible.$$
Now we estimate
$$||\phi_{Q,j} P_j u||_{L^{\infty}} \lesssim 2^{{3j \over 2}} u_{Q} +
negligible,$$
$$||\nabla \phi_{Q_j,j}\tilde  P_k u||_{L^2} \lesssim 
2^k \sum_{l=-2}^2 ||\phi_{Q_j,j} P_{k+l} u||_{L^2}
+negligible.$$
Combining these estimates gives the desired inequality.  
\end{proof} 

In fact, we get a slightly better estimates, for instance by the div-curl lemma, but
it does not seem to be necessary.

The previous three lemmas are somewhat wasteful, particularly in the definition
of the local part. We are using more squares than we need to cover. For the
result of the final section we need a somewhat more efficient decomposition,
which is proved in exactly the same way. We formulate it here.

For $Q$ a cube at level $j$, we define the ancestors of $Q$, 
$${\cal A}(Q) =\{Q_k\}_{\epsilon j <k<j-4}.$$
For $k \leq j-4$, we define the collection ${\cal S}_k(Q)$ to be
a covering (with overlap $\lesssim 1$) of the cube $Q_{j-4}$ by cubes
at level $k$. We define the strict extended family
$${\cal E(Q)} =\bigcup_{k=j-4}^{j+{1000 \over \epsilon}} {\cal S}_k(Q),$$
and we define ${\cal F(Q)}$, the followers of $Q$ by
$${\cal F(Q)}=\bigcup_{k > j + {1000 \over \epsilon}} {\cal S}_k(Q).$$

Then we have the following corollary of the proofs 
(but not the statements) of lemmas \ref{highlowlowhigh},
\ref{neighbours} and \ref{highhigh}

\begin{corollary}\label{theworks}
For any $\delta > 0$ we can estimate
$$G_Q \lesssim Z_{Q,lhhl} + Z_{Q,loc} + Z_{Q,hh} + 2^{j(1 + {3 \over 2}\delta)} u_Q
u_{{\cal N}^1(Q)},$$
where
$$Z_{Q,lhhl} = \sum_{k < j-4} \sum_{Q^{\prime} \in {\cal N}^1(Q)}
2^{j+{3k \over 2}} u_{Q_k} u_{Q^{\prime}} u_Q,$$
and
$$Z_{Q,loc} = 2^{{5j \over 2}} \sum_{Q^{\prime},Q^{\prime \prime} \in {\cal E}(Q)}
u_Q u_{Q^{\prime}} u_{Q^{\prime \prime}},$$
and finally
$$Z_{Q,hh} =\sum_{k > j + {1000 \over \epsilon}} 2^{{3j \over 2} + k} u_Q
\sum_{Q^{\prime} \in {\cal S}_k(Q)} u_{Q^{\prime}}^2.$$
\end{corollary}

Now we are ready to describe our singular set. We will say that a cube $Q$
of sidelength $2^{-j(1-2\epsilon)}$ is {\bf bad} if
\be{bad}\int_0^T \int  \sum_{k \geq j} 2^{2\alpha k} |\phi_{Q,k} P_k u|^2
 \gtrsim 2^{ -(5-4\alpha )j - 
100\epsilon j } . \end{equation}
Let $E_j$ be the union of $2^{{3000 \over \epsilon}}Q$ for all cubes  
$Q$ of sidelength $2^{-j(1-2\epsilon)}$ 
which are bad.

We will need the following well known covering lemma of Vitali (see \cite{Stein}):

\begin{lemma}\label{Vitali} Let ${\cal C}$ be any collection of cubes, then there is
a subcollection ${\cal C}^{\prime}$ so that any two cubes in 
${\cal C}^{\prime}$ are pairwise disjoint and so that
$$\bigcup_{Q \in {\cal C}} Q \subset \bigcup_{Q \in {\cal C}^{\prime}} 5Q.$$
\end{lemma} 

\begin{proposition} There is covering ${\cal Q}_j$ of $E_j$ by
cubes of sidelength $2^{-j(1-2\epsilon)}$ so that
$$\# ({\cal Q}_j) \lesssim  2^{(5-4\alpha )j + 100\epsilon j }.$$
\end{proposition}

\begin{proof}
Let ${\cal C}$ be the collection of  cubes $2Q$  where $Q$ is a
bad cube at level $j$. 
From Lemma \ref{Vitali} we know that there are disjoint cubes 
$\{ 2Q_{\rho} \}_{\rho \in \Bbb Z}$ such that the collection consisting of 
$\{ 10Q_{\rho} \}$ covers the set
$E_{j}$. Any cube of sidelength $10 \cdot 2^{-j(1-2\epsilon)}$ can be
covered by
1000 cubes of sidelength $2^{-j(1-2\epsilon)}$ (and 1000 is a constant.)
We will define ${\cal Q}_j$ to be the covering formed by the union of the
thousandths of the elements of  $\{ 10 Q_{\rho} \}$ .

In order to count cubes in ${\cal Q}_j$ it is enough to count the disjoint
cubes used in the construction of ${\cal Q}_j$.

However we know since the $2Q_\rho$'s are disjoint that
\be{gooddis} \sum_{\rho} \sum_k\int_0^T 2^{2 \alpha k} |\phi_{Q_{\rho},k} 
P_k u|^2 \lesssim
\int_0^T \int |{\Delta}^{\alpha}u|^2  \lesssim 1
\end{equation} 
by conservation of energy, while we know
\be{baddis} \sum_{\rho} \sum_{k} \int_0^T 2^{2 \alpha k}
|\phi_{Q_{\rho},k} P_j u|^2 
\gtrsim \#({\cal Q}_j) 2^{ -(5-4\alpha )j - 
100\epsilon j }  \end{equation}
by the badness of the cubes.

Combining the inequalities \eqref{gooddis} and
\eqref{baddis} implies the claim. \end{proof}

We pause for a moment to apply Lemma \ref{limsup}.

\begin{corollary} \label{dimension} The dimension of 
$E = \limsup_{j \longrightarrow \infty} E_j$
is bounded by $5-4\alpha + O(\epsilon)$. \end{corollary}

If we could show that if $x \notin E$ then $x$ is a regular point of $u$
at time $T$ - in other words that 
$$\limsup_{t \longrightarrow T} |u(x,t)| < \infty.$$
(Indeed if we could show it is a regular point for a derivative of $u$ of
any fixed order) this would immediately imply Theorem \ref{main}.

Unfortunately, we will show something slightly more complicated - the
same statement for a somewhat larger collection ${\cal Q}^{\prime}_j$.

In the following section, we will investigate what are the immediate
consequences of $x \notin E$.
Saying that $x \notin E$ is the same as saying that there exists a $j$
so that so that for any $k > j$, we have $x \notin E_k$. Denote
$F_j$ as the set of points with this property. To prove
a regularity statement about $E$, it suffices to show that that statement
holds for any $j_0$, provided $x \in F_{j_0}$.
However fixing $j$, we may change our constants in the definition
of bad square so that $x$ is not contained in any bad squares. Thus
we may as well assume $x$ is contained in no bad squares.  This will be
our hypothesis in the following section. (And as promised, the constants
will now depend on $j_0$.)

\section{Critical regularity}

Heuristically, the worst part of $G_Q$ is $G_{Q,loc}$ which looks at
scale $j$ like $2^{5j \over 2} u_Q^3$. On the other hand, dissipation
gives decay like $2^{2 \alpha j} u_Q^2$. This should mean that as long
as $u_Q < 2^{2 (\alpha - {5 \over 4}) j},$ the growth of $u_Q$ ought
to be under control. It is this estimate which we will show in this section
for any $Q$ not contained in $E_k$ for any $k > j_0$.

\begin{theorem}\label{critical} Let $Q$ be as above,
then there is a constant $C$, depending only on $T$, the
initial conditions for $u$, the constant in the definition of badness,
and $j_0$ so
that 
$$u_Q (t)<  C 2^{2 (\alpha - {5 \over 4} -{\epsilon \over 2} ) j(Q)}.$$
\end{theorem}

\begin{proof} We proceed by contradiction. Suppose the theorem
is false. We let $T_0$ be the first time and $Q$ be the largest cube so that 

$$u_Q (T_0) >  C 2^{2 (\alpha - {5 \over 4} -{\epsilon \over 2}) j}.$$

Now since our initial data is smooth, at the initial time, (since $j$ is
chosen sufficiently large), we have
$$u_Q(0) \lesssim 2^{-1000j}.$$

Therefore, it must be the case that
$$\int_0^{T_0} (G_Q(t) 
-c 2^{2\alpha j}u_Q^2 + C 2^{(2 \alpha -2\epsilon)j}
\sum_{Q^{\prime} \in {\cal N}^{1}(Q)} u_{Q^{\prime}}^2) dt \gtrsim 
2^{4 (\alpha - {5 \over 4} -{\epsilon \over 2}) j}.$$

By using proposition \ref{goodneighbours}, we can replace $Q$ by an extended
nuclear family ${\cal N}^l(Q)$, with $l < {400 \over \epsilon}$
 for which
$$\int_0^T ( -c 2^{2\alpha j}u_{{\cal N}^l(Q)}^2 + C 2^{(2 \alpha
-2\epsilon)j}
u_{{\cal N}^{l+1}(Q)}^2 )dt \lesssim - 
\int_0^T  2^{2\alpha j}u_{{\cal N}^l(Q)}^2 .$$

For the current theorem, this is all we need. Since $u_{{\cal N}^l(Q)}$
also begins $\lesssim 2^{-1000j}$, we must have

$$\int_0^{T_0} G_{{\cal N}^l(Q)}(t) dt \gtrsim 2^{4 (\alpha - {5 \over 4}
-{\epsilon \over 2} ) j}+ \int_0^T  2^{2\alpha j}u_{{\cal N}^l(Q)}^2 ,$$
where we define
$$G_{{\cal N}^l(Q)}(t) = \sum_{Q_1 \in {\cal N}^l(Q)} G_{Q_1}(t).$$
Since there are only $\lesssim 1$ cubes in ${\cal N}^l(Q)$,
for this to be the case, there must be a $\tilde Q$ so that
$$\int_0^{T_0} G_{\tilde Q}(t) dt \gtrsim 2^{4 (\alpha - {5 \over 4}
-{\epsilon \over 2} ) j}+ \int_0^T  2^{2\alpha j}
u_{{\cal N}^l(\tilde Q)}^2 ,$$

We contradict this by using Lemmas \ref{highlowlowhigh}, \ref{neighbours},
and  \ref{highhigh} to estimate $G_{\tilde Q}$. By our definition of $E_j$, the cube
$\tilde Q$ is contained in no larger bad squares.

Each estimate contains
a factor of $u_{\tilde Q}$ which we take out using the estimate
$|u_{\tilde Q}|\lesssim  2^{2(\alpha -{5 \over 4} - {\epsilon \over 2})j}$, which we
get from the definition of $T_0$.

First let $Q_1$ be a nuclear family member of $\tilde Q$ and $Q_2$ be a distant
ancestor at level $k$. 
Suppose $k < \delta j$, we see that
$2^{(1 + {3 \delta \over 2})j} u_{Q_1} \lesssim 2^{2\alpha j} u_{{\cal N}^l(\tilde Q)}$.
Thus we need not worry about this term.

Suppose $k > \delta j$.
We must estimate
$$\int_0^T 2^{{3k \over 2} + j} u_{Q_1} u_{Q_2} dt.$$ Note that
both $Q_1$ and $Q_2$ are good squares. Thus we have the estimates
$$\int_0^T 2^{3k} u_{Q_2}^2 \leq 2^{(2\alpha-2)k -10 \epsilon k}.$$
and
$$\int_0^T 2^{2j} u_{Q_1}^2 \leq 2^{(2\alpha-3)j -10 \epsilon j}.$$
Applying Cauchy-Schwartz, we get (using $\alpha \geq 1$)
$$\int_0^T 2^{{3k \over 2} + j} u_{Q_1} u_{Q_2} dt \leq  2^{({4\alpha-5
\over 2})j
-10\epsilon j}.$$
Summing over   (there are only $j$ terms) provides the desired estimate
on $G_{\tilde Q,hl} + G_{\tilde Q,lh}$. We can estimate
 $G_{\tilde Q,loc}$ in the same way (by allowing
$k$ as large as $j +{1000 \over \epsilon}$.)

We are left to estimate $G_{\tilde Q,hh}$. We fix a scale $k>j$ and pick
$k_2$ within 2 of $k$. Now we are left to estimate
$$\int_{0}^{T_0} 2^{{3j \over 2} + k}   ||\phi_{\tilde Q,k} P_k u||_{L^2}
 ||\phi_{\tilde Q,k_2} P_{k_2} u||_{L^2}.$$
However since $\tilde Q$ is a good square, 
we have
$$\int 2^{2\alpha k} |\phi_{\tilde Q,k}P_k u |^2 \leq 2^{(4 \alpha - 5)j-10 \epsilon j}.$$
Thus
$$\int 2^{k + {3j \over 2}} |\phi_{\tilde Q,k}P_k u |^2 \leq
2^{({4 \alpha- 5 \over 2})j - 10 \epsilon j - (2 \alpha-1) (k-j)},$$
which is an estimate  that decays geometrically in $k$ when $\alpha >{{1 \over 2}}$.
By using the similar estimate for $k_2$, applying Cauchy Schwartz and
summing over $k$, we get the desired result.
\end{proof}

\section{Combinatorics}\label{comb}

In this section, we begin with the set $E_j$ which is the union of a collection
${\cal Q}_j$ of cubes with sidelength $2^{-(1-2\epsilon) j}$ having cardinality
$\lesssim 2^{(5-4\alpha + 100 \epsilon) j}$. (We assume $4\alpha -4 > 200 \epsilon$.)

\begin{theorem}\label{nopaths} There exists a 
sequence of collections ${\cal Q}^{\prime}_j$ of cubes of sidelength 
$2^{-(1-2\epsilon) j}$  with $\# ( {\cal Q}^{\prime}_j) \lesssim  
2^{(5-4\alpha + 100 \epsilon) j}$, so that for any $Q$ of length
$2^{-(1-2\epsilon) j}$  which does not intersect any element 
${\cal Q}^{\prime}_j$  there exists a number ${1 \over 2} < r < 1$
with the following property:
For no $k > j$ is there $\tilde Q \in {\cal Q}_k$ so that
\be{zerostep} 100 \tilde Q \cap \partial (rQ) \neq \emptyset 
\end{equation}
\end{theorem}

\begin{proof} We refer to the elements of ${\cal Q}_j$ as the bad cubes.
We say that a cube of sidelength $2^{-(1-2\epsilon) j}$  is very bad
if either it intersects a bad cube of the same length or it intersects
more than $c 2^{(5-4\alpha + 150 \epsilon) (k-j)}$ elements of
$100{\cal Q}_k$ for some $k>j$ with $c$ a small constant to be
specified later. Let $E^{\prime}_j$ be the union of all very bad cubes
of length $2^{-(1-2\epsilon) j}$ , then by the estimates on the cardinality
of the ${\cal Q}_k$'s and by the Vitali lemma, we can see that that 
$E^{\prime}_j$ can be covered by $\lesssim 2^{(5-4\alpha + 100 \epsilon) j}$
cubes of length $2^{-(1-2\epsilon) j}$. We refer to these cubes
as ${\cal Q}^{\prime}_k$. Now we need only prove \eqref{zerostep}.

Let $Q$ be a cube of length $2^{-(1-2\epsilon) j}$  which does not intersect     
$E^{\prime}_j$. Let ${\cal D}^k(Q)$ be the set of elements of $100 {\cal Q}_k$
which intersect $Q$. Then we have the estimate
$$\#({\cal D}^k(Q)) \leq  c 2^{(5-4\alpha + 150 \epsilon) (k-j)}.$$

Let $f_k(r)$ be the function defined from ${1 \over 2}$ to $1$ which counts how many
elements of ${\cal D}^k(Q)$ intersect $\partial (r Q)$. For each 
$Q^{\prime} \in  {\cal D}^k(Q)$ define $r_{Q^{\prime}}$ to be that number so
that the center of $Q^{\prime}$ lies on $\partial (r_{Q^{\prime}} Q)$. Then
$$f_k(r) \leq \sum_{Q^{\prime} \in {\cal D}^k(Q)} 
\chi_{( r_{Q^{\prime}} - 100 (2^{(1-2\epsilon)(j-k)}), r_{Q^{\prime}}  + 
100 (2^{(1-2 \epsilon) (j-k)}) )    }.$$
Thus
$$||f_k(r)||_{L^1} \leq 200 c 2^{ (4\alpha - 4 - 152 \epsilon) (j-k)}.$$
Since $4\alpha - 4 > 200 \epsilon,$ this estimate decays geometrically
with $k$. By choosing $c$ sufficiently small, we may arrange that
$$|| \sum_{k>j} f_k||_{L^1} < {1 \over 4}.$$
Thus by Tchebychev's inequality, there must be a value of $r$ between
${1 \over 2}$ and $1$ so that $f_k(r)=0$ for all $k > j$. This
is the value of $r$ that we choose.
\end{proof} 

\section{Barrier estimate}

In this section, we prove regularity on the interior of a cube $Q$, provided
that one has critical regularity for cubes containing it and cubes $Q^{\prime}$ for
which $\partial Q \cap 100 Q^{\prime} \neq \emptyset$.

If $Q$ is a cube  and $Q_{1} \subset Q$,
we define $d(Q_{1})$, the graph distance of $Q_{1}$ to the boundary
of $Q$ by $d(Q_{1})=k-1$, where $k$ is the smallest positive integer so
that
$$2^k Q_{1} \cap \partial Q \neq \emptyset.$$

\begin{lemma} \label{Barrier-estimate} Let $Q$ be a cube. Suppose we know that
for all $t < T$ we have that for any cube $Q^{\prime}$ so that 
$Q \subset Q^{\prime}$
with sidelength of $Q^{\prime}$ being $2^{-l(1-2\epsilon)}$, we have that
$$|u_{Q^{\prime}}(t)| \lesssim 2^{-l ({5-4 \alpha + 2\epsilon \over
2})},$$
and suppose further that for any $Q^{\prime}$ so that 
$100 Q^{\prime} \cap \partial Q \neq \emptyset$
with sidelength of $Q^{\prime}$ being $2^{-l(1-2\epsilon)}$
and $l > j-2$, we have that
$$|u_{Q^{\prime}}(t)| \lesssim 2^{-l ({5-4 \alpha + 2\epsilon \over 2})},$$
then for any $Q_1 \subset Q$ of length $2^{-k(1-2\epsilon)}$, we have the estimate
\be{induction}|u_{Q_1}(t)| \lesssim 2^{-k \rho(Q_1)}, \end{equation}
where
$$\rho(Q_1) = \min (10, {5-4 \alpha + 2 \epsilon \over 2} + 
{\epsilon (d(Q_1)-5) \over 50}).$$
\end{lemma}

\begin{proof} At time $t=0$ for all $Q_1$ we have
$$u_{Q_1}(t) \lesssim 2^{-1000 k},$$
when $Q_1$ has sidelength $2^{-(1-2\epsilon)k}$. Thus at time 0, the lemma
is satisfied. Let $t_0$ be the first time at which the lemma fails and $Q_1$
be one of the cubes for which it fails. It must be the case by hypothesis
that $32 Q_1$ does not intersect $\partial Q$.
 Then we have
$$u_{Q_1}(t_0) \sim 2^{-k \rho(Q_1)}.$$
Let $t_1 < t_0,$ be the last time before $t_0$ when
$$u_{Q_1}(t_1) \lesssim 2^{k (-\rho(Q_1) +{\epsilon \over 10})}.$$
Then we have
\be{goalline} \int_{t_1}^{t_0} {d \over dt} (u_{Q_1}^2) \gtrsim 
2^{-2k \rho(Q_1)}.
\end{equation}
However on the time interval $(t_1,t_0)$ the lemma is satisfied.

We will invoke Corollary \ref{theworks}.
Now for any $Q_2 \in {\cal E}(Q_1)$ we have 
\be{ext} d(Q_2) \geq d(Q_1) - 5.\end{equation}
This is because $Q_2 \subset 32 Q_1$.
Further for any ancestor $Q_3 \in {\cal A}(Q_1)$ with $Q_3$ having sidelength $2^{-(1-2\epsilon)l}$, with $\epsilon k < l < k-4$, 
\be{anc}d(Q_3) \geq d(Q_1) - 5 (k-l). \end{equation}
Further for any follower  $Q_4 \in {\cal F}(Q_1)$ (that is a cube which 
contributes to $G_{Q_1,hh}$ and is in particular contained in 
${3 \over 2} Q_1$)
of $Q_1$ with $Q_4$ having sidelength
$2^{-(1-2\epsilon)l}$ with $l >k+ {1000 \over \epsilon}$, we have
 \be{dec}d(Q_4) \geq d(Q_1) +  ({l-k \over 2}). \end{equation}
Applying \eqref{dec}, we conclude that $\rho(Q_4)=10$.
Thus we conclude

$$Z_{Q_1,hh} \lesssim \sum_{l>k} 2^{{5k \over 2}} 2^{3(l-k)} 2^{-20 l} 
2^{(-\rho(Q_1) + {\epsilon \over 10})k},$$
by counting the elements of ${\cal S}_l(Q)$ by  $2^{3(l-k)}$. Calculating, we 
find
$$Z_{Q_1,hh} \lesssim 2^{-({35 \over 2} + \rho(Q_1)) k},$$
which since $\rho(Q_1) \leq 10,$ cannot possibly account for \eqref{goalline}.

Now applying \eqref{ext} to Proposition \ref{diss}, we observe that for any $t$
between $t_1$ and $t_0$, we have dissipation 
at $Q_1$ of $\gtrsim 2^{2 \alpha k} u_{Q_1}^2(t) \gtrsim  
2^{(2 \alpha - 2 \rho(Q_1) + {\epsilon \over 5})k}$.
Thus to reach a contradiction, it suffices to show that on this
time interval
$$Z_{Q_1,hllh} + Z_{Q_1,loc} \lesssim   
2^{(2 \alpha - 2 \rho(Q_1) + {\epsilon \over 5})k}.$$
For this reason we can also ignore the ultra-low term 
$2^{j (1 + {3 \over 2} \delta)} u_{Q_1} u_{{\cal N}^1(Q_1)}$
in Corollary \ref{theworks}.

Using \eqref{ext}, we observe that for any $Q_2$ in  ${\cal E}(Q_1)$
we have the estimates
$$u_{Q_2} \lesssim 2^{(-\rho(Q_1) +{\epsilon \over 10})k},$$
as well as
$$u_{Q_2} \lesssim 2^{ {(4\alpha -5 - 2\epsilon)k  \over 2}},$$
by the lower bound on $\rho$.
Now we calculate
$$Z_{Q_1,loc} \lesssim \sum_{Q_2,Q_2^{\prime} \in {\cal E} (Q_1)} 2^{{5k \over 2}}
u_{Q_2} u_{Q_2^{\prime}} u_{Q_1}
 \lesssim 2^{(2 \alpha -{4 \epsilon \over 5}-2 \rho(Q_1)) k}.$$
Thus $Z_{Q_1,loc}$ cannot contribute to the growth.

Now to estimate $Z_{Q_1,hllh}$, we observe that for $Q_1^{\prime}
\in {\cal N}^1(Q_1)$ we have the estimate
$$u_{Q_1^{\prime}} \lesssim 2^{(-\rho(Q_1) +{\epsilon \over 10})k},$$
while for ancestor $Q_4$ of sidelength $2^{-(1-2\epsilon)l}$, 
we apply \eqref{anc} (as well as the hypotheses of the lemma for squares 
larger than $Q$)
$$u_{Q_4} \lesssim 2^{ {(4\alpha -5 - 2\epsilon)k  \over 2}+
 ({3 \over 2}-100 \epsilon)(k-l)},$$
since of course ${3 \over 2} > {5-4\alpha  \over 2} + 200 \epsilon$.
Now we just estimate
$$Z_{Q_1,hllh} \lesssim \sum_l \sum 2^{{3l \over 2} + k} u_{Q_1^{\prime}}  
u_{Q_4} u_{Q_1} \lesssim 2^{(2 \alpha -{4 \epsilon \over 5}-2\rho(Q_1)) k}.$$
Thus $u_{Q_1}$ could not have grown which is a contradiction. \end{proof}

But now we we have in fact proven the main theorem. We define $E_j^{\prime}$ as in
section \ref{comb}. We need to show that $u$ is regular at any $x$ not
contained in any $E_j^{\prime}$  with $j$ larger than some integer $j_x$. By changing 
our constants, we can say $x$ is not in in any $E_j^{\prime}$. Now by
theorems \ref{critical} and \ref{nopaths}, for any cube $Q$ centered at $x$,
we can find a cube almost half as large which satisfies the hypotheses of
Lemma \ref{Barrier-estimate}. But the conclusion of the lemma implies that 
$x$ is a regular point. Thus any singular point must be contained in
$E=\limsup E_j^{\prime}$. By lemma \ref{limsup}, we have $\dim(E) < 5 - 4\alpha +
20 \epsilon$. Letting $\epsilon$ tend towards 0, we get Theorem \ref{main}.

\end{document}